\documentclass[12pt]{article}
\usepackage{amsfonts}
\usepackage{amssymb}
\usepackage{amsthm}
\usepackage{enumerate}
\usepackage{amsmath}
\usepackage{times}
\usepackage{amscd}
\usepackage[all]{xy}
\usepackage{verbatim}
\usepackage{graphicx}
\usepackage{isolatin1}
\usepackage{hyperref}
\usepackage{longtable}
\title{On the computation of Galois representations associated to level one modular forms}
\author{Johan Bosman
\footnote{This research was partially supported by the Dutch scientific organisation NWO.
\ \ \ \ \ \ \ \ \ \ \ \ \ \ \ \ \ \ \ \ \ \ \ \
E-mail: \mbox{jgbosman@math.leidenuniv.nl}}}
\date{}

\theoremstyle{plain}
\newtheorem{theorem}{Theorem}
\newtheorem{lemma}{Lemma}
\newtheorem{corollary}{Corollary}
\newtheorem{proposition}{Proposition}

\theoremstyle{definition}

\newenvironment{proof2}[1]{\medskip \noindent {\bf #1}}
           {\unskip \nobreak \hfill \hbox{$\Box$} \par \medskip}

\newtheorem{exam}{Voorbeeld}

\renewcommand\mod{\mathop{\rm mod} \nolimits}
\newcommand\Frob{\mathop{\rm Frob} \nolimits}

\newcommand\F{\mathbb{F}}
\newcommand\Z{\mathbb{Z}}
\newcommand\Q{\mathbb{Q}}

\newcommand\C{\mathbb{C}}

\newcommand\Qbar{\overline{\mathbb{Q}}}
\renewcommand\P{\mathbb{P}}

\newcommand\el{\ell}

\newcommand\jac[2]{\left(\frac{#1}{#2}\right)}

\newcommand\GL{\mathop{\rm GL} \nolimits}
\newcommand\SL{\mathop{\rm SL} \nolimits}
\newcommand\PGL{\mathop{\rm PGL} \nolimits}

\newcommand\Gal{\mathop{\rm Gal} \nolimits}

\renewcommand\ker{\mathop{\rm Ker} \nolimits}
\newcommand\im{\mathop{\rm Im} \nolimits}

\newcommand\tbox[1]{\mathop{\rm #1} \nolimits}
\newcommand\surj\twoheadrightarrow

\newcommand\disc{\mathop{\rm Disc} \nolimits}
\newcommand\tr{\mathop{\rm tr} \nolimits}
\renewcommand\O{\mathcal{O}}
\renewcommand\o[1]{\overline{#1}}

\newcommand{\GQl}{\Gal(\o\Q_\el/\Q_\el)}
\newcommand{\GQ}{\Gal(\Qbar/\Q)}
\newcommand{\GQp}{\Gal(\o\Q_p/\Q_p)}
\newcommand{\m}{{{}^{\underline{\phantom{\circ}}}}}
\newcommand{\mm}{{{}^{\underline{\phantom{\cdot}}}}}
\newcommand\mat[4]{
  \left(
  {#1 \atop #3}
  \thinspace
  {#2 \atop #4}
  \right)
}

\makeatletter
\def\eqalign#1{\null\,\vcenter{\openup\jot\m@th
  \ialign{\strut\hfil$\displaystyle{##}$&$\displaystyle{{}##}$\hfil
      \crcr#1\crcr}}\,}

\newcounter{Lcount}
{%
  \medskip
  \setcounter{Lcount}{1}%
  \begin{list}{\bf\arabic{section}.\arabic{Lcount}.}{\usecounter{Lcount}%
    \settowidth\leftmargin{\bf\arabic{section}.99. }%
    \setlength\rightmargin{0pt}%
    \setlength\itemsep{\bigskipamount}}%
}%
{\end{list}}

{%
  \par
  \setlength{\leftmargin}{#1}%
  \setlength{\rightmargin}{#2}%
  \advance\linewidth -\leftmargin
  \advance\linewidth -\rightmargin
  \advance\@totalleftmargin\leftmargin
  \@setpar{{\@@par}}%
  \parshape 1 \@totalleftmargin \linewidth
  \ignorespaces
}%
{\par}

\makeatother
\begin{document}
\maketitle
\begin{abstract}
\noindent
In this paper we 
explicitly compute mod-$\el$ Galois representations associated to modular forms. To be precise,
we look at cases with $\el\leq 23$ and the modular forms considered will be cusp forms of level $1$ and weight up to~$22$.
We present the result in terms of polynomials associated to the projectivised representations.
As an application, we will improve a known result on Lehmer's non-vanishing conjecture for Ramanujan's tau function.
\end{abstract}
\section{Introduction}
The Ramanujan tau function is the function $\tau:\Z_{>0}\to\Z$ defined by
$$
\Delta=q\prod_{n\geq 1}\left(1-q^n\right)^{24} = \sum_{n\geq 1}\tau(n)q^n.
$$
If we write $q=\exp(2\pi iz)$ for $z$ in the complex upper half plane 
then $\Delta(z)$ is a holomorphic cusp form of level $1$ and weight $12$.
We have the relations 
$$
\begin{array}{rcll}
\tau(mn) &=& \tau(m)\tau(n) &\quad\mbox{if $\gcd(m,n)=1$},\cr
\tau(p^{r+1}) &=& \tau(p)\tau(p^r)-p^{11}\tau(p^{r-1})&\quad\mbox{for $p$ prime and $r\geq 1$}.
\end{array}
$$
These relations determine $\tau(n)$ in terms of $\tau(p)$ for $p$ prime. 
\newline

For $\el\in\{2,3,5,7,23,691\}$ there exist simple formulas for $\tau(p)\mod\el$, or in some cases even modulo certain powers of $\el$;
for instance $\tau(p)\equiv p^{41}+p^{70}\mod 5^3$ for primes $p\not=5$ and $\tau(p)\equiv 1+p^{11}\mod 691$ for all primes $p$. 
In general, there is a Galois representation $\rho=\rho_{\Delta,\el}:\Gal(\Qbar/\Q)\to\GL_2(\F_\el)$ unramified outside $\el$
such that for any Frobenius element $\Frob_p\in\Gal(\Qbar/\Q)$ attached to a prime $p\not=\el$ the characteristic polynomial
of $\rho(\Frob_p)$ is congruent to \mbox{$X^2-\tau(p)X+p^{11}$} modulo $\el$.
The simple congruences for special values of $\el$ are due to the
fact that the image of $\rho$ does not contain $\SL_2(\F_\el)$ in those cases; 
such a representation is called \emph{exceptional} and is in many cases easy to compute. 
\newline

Besides the modular form $\Delta$ of weight $12$ 
we will also consider the unique normalised cusp forms of level $1$ and weights $16$, $18$, $20$ and $22$ in this paper. 
To fix a notation, for any $k\in\Z$
satisfying $\dim S_k(\Gamma(1))=1$ we will denote the unique normalised cusp form in $S_k(\Gamma(1))$ by $\Delta_k$. 
We will denote the coefficients of the $q$-expansion of $\Delta_k$ by $\tau_k(n)$:
$$
\Delta_k(z)=\sum_{n\geq 1} \tau_k(n)q^n\in S_k(\Gamma(1)).
$$
From $\dim S_k(\Gamma(1))=1$ it follows that the numbers $\tau_k(n)$ are integers. 
For every $\Delta_k$ 
and every prime $\el$ there is a representation \mbox{$\rho_{\Delta_k,\el}:\GQ\to\GL_2(\F_\el)$} such that for every prime $p\not=\el$ 
we have that the characteristic polynomial of $\rho_{\Delta_k,\el}(\Frob_p)$ is congruent to \mbox{$X^2-\tau_k(p)X+p^{k-1}\mod\el$}.
For a summary on the exceptional representations $\rho_{\Delta_k,\el}$ and the corresponding congruences for $\tau_k(n)$, see~\cite{SwD}.
\newline

In this paper we shall present polynomials that belong to the projectivisations of the non-exceptional
Galois representations belonging to rational level
one forms modulo primes up to $23$. 
Finding these polynomials is a matter of experimental computation, but the known cases of Serre's conjecture 
permit us to verify the correctness.
As a by-product we will verify Lehmer's conjecture of the non-vanishing of $\tau(n)$ (see \cite[p. 429]{Le}) to a higher bound than 
what was done before.

\subsection{Notational conventions}
Throughout this paper, for every field $K$ we will fix an algebraic closure $\o K$ and all algebraic extension fields of $K$ will be regarded as subfields of 
$\o K$. Furthermore, for each prime number $p$ we will fix an embedding $\Qbar\hookrightarrow\Qbar_p$ and hence an embedding
$\GQp\hookrightarrow\GQ$, whose image we call $D_p$. 
We will use $I_p$ to denote the inertia subgroup of $\GQp$.  
\newline

All representations (either linear or projective) in this paper will be \emph{continuous}.
For any field $K$, a linear representation 
$\rho:G\to\GL_n(K)$ defines a 
projective representation
$\tilde\rho:G\to\PGL_n(K)$ via the canonical map $\GL_n(K)\to\PGL_n(K)$. 
We say that a projective representation $\tilde\rho:G\to\PGL_n(K)$ is \emph{irreducible} if the induced action of $G$ on $\P^{n-1}(K)$ fixes
no proper subspace. So for $n=2$ this means that every point of $\P^1(K)$ has its stabiliser subgroup not equal to $G$.

\subsection{Statement of results}
\begin{proposition}\label{PGLprop}
For every pair $(k,\el)$ occurring in Table \ref{polytable} on page \pageref{polytable}, let the polynomial $P_{k,\el}$ be defined as in 
that same table.
Then the splitting field of each $P_{k,\el}$ is the fixed field of $\ker(\tilde\rho_{\Delta_k,\el})$ and has Galois group $\PGL_2(\F_\el)$.
Furthermore, if $\alpha\in\Qbar$ is a root of $P_{k,\el}$ then the subgroup of $\Gal(\Qbar/\Q)$ fixing $\alpha$ corresponds via 
$\tilde\rho_{\Delta_k,\el}$ to a subgroup of $\PGL_2(\F_\el)$ fixing a point of $\P^1(\F_\el)$.
\end{proposition}
For completeness we also included the pairs $(k,\el)$ for which $\rho_{k,\el}$ is isomorphic to the action of $\GQ$ on the $\el$-torsion of an elliptic curve.
These are the pairs in Table \ref{polytable} with $\el=k-1$, as there the representation is the $\el$-torsion of $J_0(\el)$, which happens to be an 
elliptic curve for $\el\in\{11,17,19\}$.
A simple calculation with division polynomials \cite[Chapter II]{LaEl} can be used to treat these cases. 
In the general case, one has to work in the more complicated Jacobian variety $J_1(\el)$, which has dimension $12$ for $\el=23$ for instance.
\newline

We can apply Proposition \ref{PGLprop} to verify the following result.
\begin{corollary}\label{Lehmer}
The non-vanishing of $\tau(n)$ holds for all 
$$
n<22798241520242687999\approx 2\cdot 10^{19}.
$$
\end{corollary}
The non-vanishing of $\tau(n)$ was verified for all $n<22689242781695999\approx 2\cdot 10^{16}$ in \cite{JoKe}.
\newline

To compute the polynomials, the author used a weakened version of algorithms described in \mbox{\cite[Sections 11 \& 24]{ECJ}}. 
The used algorithms do not give a proven output,
so we have to concentrate on the verification.
We will show how to verify the correctness of the polynomials in Section \ref{secproof} after setting up some preliminaries about Galois representations
in Section \ref{secgalrep}. In Section \ref{secLehmer} we will point out how to use Proposition~\ref{PGLprop} in a calculation that verifies Corollary
\ref{Lehmer}. All the calculations were perfomed using {\sc Magma} (see \cite{Magma}).
 
\section{Galois representations}\label{secgalrep}
This section will be used to state some results on Galois representations that we will need in the proof of Proposition \ref{PGLprop}.
\subsection{Liftings of projective representations}
Let $G$ be a topological group, let $K$ be a field and let $\tilde\rho:G\to\PGL_n(K)$ be a projective representation. Let $L$ be an extension field of $K$.
By a \emph{lifting} of $\tilde\rho$ over $L$ we shall mean a representation $\rho:G\to\GL_n(L)$ that makes the following diagram commute:
$$
\SelectTips{cm}{10}
\xymatrix{
G \ar [r]^-{\tilde\rho} \ar [d]_{\rho} & \PGL_n(K) \ar @{^(->}[d]\cr
\GL_n(L) \ar @{->>}[r] & \PGL_n(L)
}
$$
where the maps on the bottom and the right are the canonical ones. If the field $L$ is not specified then by a lifting of $\tilde\rho$ we shall mean a lifting
over $\o K$.
\newline 

An important theorem of Tate arises in the context of liftings. For the proof we refer to \cite[Section 6]{SeWt1}. 
Note that in the reference representations over $\C$ are
considered, but the proof works for representations over arbitrary algebraically closed fields. 

\begin{theorem}[Tate]\label{tate2}
Let $K$ be a field and let $\tilde\rho:\Gal(\Qbar/\Q)\to\PGL_n(K)$ be a projective representation.
For each prime number $p$, let $\rho'_p:D_p\to\GL_n(\o K)$ be a lifting of $\tilde\rho|_{D_p}$. Assume that all but finitely many of the 
$\rho'_p$ are unramified. Then there is a unique lifting $\rho:\Gal(\Qbar/\Q)\to\GL_n(\o K)$ such that for all primes $p$ we have
$$
\rho|_{I_p}=\rho'_p|_{I_p}.
$$
\end{theorem}

\begin{lemma}\label{localunramified}
Let $p$ be a prime number and let $K$ be a field. Suppose that we are given an unramified projective representation
$\tilde\rho_p:\Gal(\o\Q_p/\Q_p)\to\PGL_n(K)$.
Then there exists a lifting 
$\rho_p:\Gal(\o\Q_p/\Q_p)\to\GL_n(\o K)$
of $\tilde\rho_p$ that is unramified as well.
\end{lemma}
\begin{proof}
An unramified homomorphism from $\Gal(\o\Q_p/\Q_p)$ to any group factors through
$\Gal(\o\F_p/\F_p)\cong\hat\Z$ and is determined by the image of
$\Frob_p\in\Gal(\o\F_p/\F_p)$. This image is an element of $\PGL_n(K)$ of finite order, say of order $m$.  
If we take any lift $F$ of $\tilde\rho(\Frob_p)$ to $\GL_n(K)$ then we have $F^m=a$ for some $a\in K^*$.
So $F':=\alpha^{\m1}F$, where $\alpha\in\o K$ is any $m$-th root of $a$, has order $m$ in $\GL_n(\o K)$.
Hence the homomorphism $\GQp\to\GL_n(\o K)$ obtained by the composition
$$
\SelectTips{cm}{10}
\entrymodifiers={+!!<0pt,\fontdimen22\textfont2>}
\xymatrixcolsep{.25pc}
\xymatrix{
\GQp \ar @{->>}[rr]
&&
\Gal(\o\F_p/\F_p) \ar [rrr]^-\sim
&&&
\hat\Z \ar @{->>}[rr]
&&
\Z/m\Z \ar [rrr]^-{1\mapsto F'}
&&&
\GL_n(\o K)
}
$$
lifts $\tilde\rho$ and is continuous as well as unramified.
\end{proof}

\subsection{Serre invariants and Serre's conjecture}
Let $\el$ be a prime. A Galois representation $\rho:\GQ\to\GL_2(\o\F_\el)$ has a \emph{level} $N(\rho)$ and a \emph{weight} $k(\rho)$. 
The definitions were introduced by Serre (see \cite[Sections 1.2 \& 2]{Se}). Later on, Edixhoven found an improved definition for the weight, which is the one we will use, 
see \cite[Section 4]{Ed}.
The level $N(\rho)$ is defined as the prime-to-$\el$ part of the Artin conductor of $\rho$ and equals $1$ if $\rho$ is unramified outside~$\el$. 
The weight is defined in terms of the local representation $\rho|_{D_\el}$; its definition is rather lenghty so we will not write it out here. When we
need results about the weight we will just state them. 
Let us for now mention that one can consider the weights of the twists $\rho\otimes\chi$ of a representation $\rho:\GQ\to\GL_2(\o\F_\el)$ by a character
$\chi:\GQ\to\o\F_\el^*$. If one chooses $\chi$ so that $k(\rho\otimes\chi)$ is minimal, then we always have $1\leq k(\rho\otimes\chi)\leq\el+1$ and 
we can in fact choose our $\chi$ to be a power of the $\mod \el$ cyclotomic character.
\newline

Serre conjectured \cite[Conjecture 3.2.4]{Se} that if $\rho$ is irreducible and odd, then $\rho$ belongs to a
modular form of level $N(\rho)$ and weight $k(\rho)$. Oddness here means that the image of a complex conjugation has determinant $\m1$. 
A proof of this conjecture in the case $N(\rho)=1$ has been published by Khare and Wintenberger:
\begin{theorem}[Khare \& Wintenberger, {\cite[Theorem 1.1]{KhWi1}}]\label{serreconj}
Let $\el$ be a prime number and let \mbox{$\rho:\Gal(\Qbar/\Q)\to\GL_2(\o\F_\el)$}
be an odd irreducible representation of level \mbox{$N(\rho)=1$}. Then there exists a modular form $f$ of level $1$ and weight
$k(\rho)$ which is a normalised eigenform and a prime $\lambda\mid\el$ of $K_f$
such that $\rho$ and $\rho_{f,\lambda}$ become isomorphic after a suitable embedding of $\F_\lambda$ into $\o\F_\el$.
\end{theorem}

\subsection{Weights and discriminants}
If a representation $\rho:\GQ\to\GL_2(\o\F_\el)$ is wildly ramified at $\el$ it is possible to relate the weight to discriminants of certain number fields.
In this subsection we will present a theorem of Moon and Taguchi on this matter and derive some results from it that are of use to us.
\begin{theorem}[Moon \& Taguchi, {\cite[Theorem 3]{MoTa}}]\label{MoonTaguchi}
Consider a wildly ramified representation $\rho:\GQl\to\GL_2(\o\F_\el)$. Let $\alpha\in\Z$ be such that 
$k(\rho\otimes\chi_\el^{\m\alpha})$ is minimal where \mbox{$\chi_\el:\GQl\to\F_\el^*$} is the mod $\el$ cyclotomic character. 
Put $\tilde{k}=k(\rho\otimes\chi_\el^{\m\alpha})$, put $d=\gcd(\alpha,\tilde{k}-1,\el-1)$ and define $m\in\Z$ by letting $\el^m$ 
be the wild ramification degree of 
\mbox{$K:=\o\Q_\el^{\ker(\rho)}$} over $\Q_\el$. Then we have
$$
v_\el(\mathcal{D}_{K/\Q_\el})=
\left\{
\begin{array}{ll}
1+\frac{\tilde{k}-1}{\el-1}-\frac{\tilde{k}-1+d}{(\el-1)\el^m}&\mbox{if $2\leq\tilde{k}\leq\el$},\cr
2+\frac{1}{(\el-1)\el}-\frac{2}{(\el-1)\el^m}&\mbox{if $\tilde{k}=\el+1$},
\end{array}
\right.
$$
where $\mathcal{D}_{K/\Q_\el}$ denotes the different of $K$ over $\Q_\el$ and $v_\el$ is normalised by $v_\el(\el)=1$.
\end{theorem}
We can simplify this formula to one which is useful in our case:
\begin{corollary}\label{cordisc1}
Let $\tilde\rho:\Gal(\Qbar/\Q)\to\PGL_2(\F_\el)$ be an irreducible projective representation that is wildly ramified at $\el$. 
Take a point in $\P^1(\F_\el)$, let $H\subset\PGL_2(\F_\el)$ be its stabiliser subgroup and let $K$ be the number field defined as
$$
K = \Qbar^{\tilde\rho^{\,\mm1}(H)}.
$$
Then the $\el$-primary part of $\disc(K/\Q)$ is related to the minimal weight $k$ of the liftings of $\tilde\rho$ by the following formula:
$$
v_\el(\disc(K/\Q))=k+\el-2.
$$
\end{corollary}
\begin{proof}
Let $\rho$ be a lifting of $\tilde\rho$ of minimal weight. 
Since $\rho$ is wildly ramified, after a suitable conjugation in $\GL_2(\o\F_\el)$ we may assume
\begin{equation}\label{rhoIl}
\rho|_{I_\el}=\mat{\chi_\el^{k-1}}{*}{0}{1},
\end{equation}
where $\chi_\el:I_\el\to\F_\el^*$ denotes the mod $\el$ cyclotomic character; this follows from the definition of weight. 
The canonical map 
$\GL_2(\o\F_\el)\to\PGL_2(\o\F_\el)$ is injective on the subgroup $\mat{*}{*}{0}{1}$, so the subfields of $\Qbar_\el$ cut out by 
$\rho$ and $\tilde\rho$ are equal, call them $K_2$. 
Also, let $K_1\subset K_2$ be the fixed field of the diagonal matrices in $\im\rho|_{I_\el}$. 
We see from (\ref{rhoIl}) that in the notation of Theorem \ref{MoonTaguchi} we can put
$\alpha=0$, $m=1$ and $d=\gcd(\el-1,k-1)$. 
So we have the following diagram of field extensions:
$$
\xymatrix{
& K_2 \ar @{-}[ld]_{\displaystyle\chi_\el^{k-1}} \ar @{-}[dd]^{\displaystyle\mat{\chi_\el^{k-1}}{*}{0}{1}}\cr
K_1 \ar @{-}[rd]_{\displaystyle\deg=\el}\cr
& \Q_\el^{\text{un}}
}
$$
The extension $K_2/K_1$ is tamely ramified of degree $(\el-1)/d$ hence we have 
$$
v_\el(\mathcal{D}_{K_2/K_1})
=\frac{(\el-1)/d-1}{(\el-1)\el/d} 
=\frac{\el-1-d}{(\el-1)\el}.
$$
Consulting Theorem \ref{MoonTaguchi} for the case $2\leq k \leq \el$ now yields
$$
\eqalign{
v_\el(\mathcal{D}_{K_1/\Q_\el^\text{un}})
&=
v_\el(\mathcal{D}_{K_2/\Q_\el^\text{un}})
-
v_\el(\mathcal{D}_{K_2/K_1})
\cr&=
1+\frac{k-1}{\el-1}-\frac{k-1+d}{(\el-1)\el}
-
\frac{\el-1-d}{(\el-1)\el}
=
\frac{k+\el-2}{\el}
}
$$
and also in the case $k=\el+1$ we get
$$
v_\el(\mathcal{D}_{K_1/\Q_\el^\text{un}})
=
2+\frac{1}{(\el-1)\el}-\frac{2}{(\el-1)\el}
-
\frac{\el-2}{(\el-1)\el}
=
\frac{k+\el-2}{\el}.
$$

Let $L$ be the number field 
$\Qbar^{\ker(\tilde\rho)}$. From the irreducibility of $\tilde\rho$ and the fact that $\im\tilde\rho$ has an element of order $\el$
it follows that the induced action of $\GQ$ on $\P^1(\F_\el)$ is transitive and
hence that $L$ is the normal closure of $K$ in $\Qbar$. This in particular implies that $K/\Q$ is wildly ramified. 
Now from $[K:\Q]=\el+1$ it follows that there are two primes in $K$ above $\el$: 
one is unramified and the other has inertia degree $1$ and ramification degree $\el$.
From the considerations above it now follows 
that any ramification subgroup of $\Gal(L/\Q)$ at $\el$ is isomorphic to a subgroup of
$\mat{*}{*}{0}{1}\subset\GL_2(\o\F_\el)$ of order $(\el-1)\el/d$ with $d\mid \el-1$. The only subgroup of index $\el$ is the subgroup of diagonal
matrices. Hence we have
$$
v_\el(\disc(K/\Q))=v_\el(\disc(K_1/\Q_\el^\text{un}))=\el\cdot v_\el(\mathcal{D}_{K_1/\Q_\el^\text{un}})=k+\el-2.
$$
\end{proof}
\begin{corollary}\label{cordisc2}
Let $\tilde\rho:\Gal(\Qbar/\Q)\to\PGL_2(\F_\el)$ be an irreducible projective representation 
and let $\rho$ be a lifting of $\tilde\rho$ of minimal weight.
Let $K$ be the number field belonging to a point of $\P^1(\F_\el)$, as in Corollary \ref{cordisc1}. 
If $k\geq 3$ is such that
$$
v_\el(\disc(K/\Q))=k+\el-2
$$
holds, then we have $k(\rho)=k$. 
\end{corollary}
\begin{proof}
From $v_\el(\disc(K/\Q))=k+\el-2\geq\el+1$ it follows that $\tilde\rho$ is wildly ramified at $\el$ so we can apply Corollary \ref{cordisc1}.
\end{proof}

\section{Proof of Proposition \ref{PGLprop}}\label{secproof}
To prove Proposition \ref{PGLprop} we need to do several verifications. We will 
derive representations from the polynomials $P_{k,\el}$ and verify that they satisfy the conditions of Theorem \ref{serreconj}. Then we know there are
modular forms attached to them that have the right level and weight and uniqueness follows then easily.
\newline

First we we will verify that the polynomials $P_{k,\el}$ from Table
\ref{polytable} have the right Galois group. The algorithm described in \cite[Algorithm 6.1]{GeKl} can be used perfectly to do this verification; proving
$A_{\el+1}\not<\Gal(P_{k,\el})$ is the most time-consuming part of the calculation here.
It turns out that in all cases we have
\begin{equation}\label{rightgalois}
\Gal(P_{k,\el})\cong \PGL_2(\F_\el).
\end{equation}
That the action of $\Gal(P_{k,\el})$ on the roots of $P_{k,\el}$ is compatible with the action of $\PGL_2(\F_\el)$ follows from the following well-known
lemma:
\begin{lemma}
Let $\el$ be a prime and let $G$ be a subgroup of $\PGL_2(\F_{\el})$ of index \mbox{$\el+1$}. 
Then $G$ is the stabiliser subgroup of a point in $\P^1(\F_\el)$. 
In particular any transitive permutation representation of $\PGL_2(\F_\el)$ of degree $\el+1$ is isomorphic to the standard action on $\P^1(\F_\el)$.
\end{lemma}
\begin{proof}
This follows from \cite[Proof of Theorem 6.25]{Su}.
\end{proof}
So now we have shown that the second assertion in Proposition \ref{PGLprop} follows from the first one. 
\newline

Next we will verify that we can obtain representations from this that have the right Serre invariants. 
Let us first note that the group $\PGL_2(\F_\el)$ has no outer automorphisms. 
This implies that for every $P_{k,\el}$, two isomorphisms (\ref{rightgalois}) define isomorphic
representations $\GQ\to\PGL_2(\F_\el)$ via composition with the canonical map $\GQ\surj\Gal(P_{k,\el})$. In other words, every $P_{k,\el}$ gives a
projective representation $\tilde\rho:\GQ\to\PGL_2(\F_\el)$ that is well-defined up to isomorphism.
\newline

Now, for each $(k,\el)$ in Table \ref{polytable}, the polynomial $P_{k,\el}$ is irreducible and hence defines a number field
$$
K_{k,\el}:=\Q[x]/(P_{k,\el}),
$$
whose ring of integers we will denote by $\O_{k,\el}$. 
It is possible to compute $\O_{k,\el}$ using the algorithm from \cite[Section 6]{BuLe} (see also \cite[Theorems 1.1 \& 1.4]{BuLe}), since we know what kind of ramification behaviour to expect.
In all cases it turns out that we have
$$
\disc(K_{k,\el}/\Q)=(-1)^{(\el-1)/2}\el^{k+\el-2}.
$$
We see that for each $(k,\el)$ the representation $\tilde\rho_{k,\el}$ is unramified outside $\el$. From Lemma \ref{localunramified} it follows that
locally outside $\el$ every $\tilde\rho_{k,\el}$ has an unramified lifting. 
Above we saw that via $\tilde\rho_{k,\el}$ the action of $\GQ$ on the set of roots of $P_{k,\el}$ is compatible with the action of $\PGL_2(\F_\el)$ on 
$\P^1(\F_\el)$, hence we can apply Corollary \ref{cordisc2} to show that the minimal weight of a lifting of $\tilde\rho_{k,\el}$ equals $k$.
Theorem \ref{tate2} now shows that every $\tilde\rho_{k,\el}$ 
has a lifting $\rho_{k,\el}$ that has level $1$ and weight $k$. From $\im\tilde\rho_{k,\el}=\PGL_2(\F_\el)$ it follows that each $\rho_{k,\el}$ is absolutely
irreducible. 
\newline

To apply Theorem \ref{serreconj} we should still verify that $\rho_{k,\el}$ is odd. Let $(k,\el)$ be given and suppose
$\rho_{k,\el}$ is even.
Then a complex conjugation $\GQ$ is sent to a matrix $M\in\GL_2(\o\F_\el)$ of determinant $1$ and of order $2$. 
Because $\el$ is odd, this means 
$M=\pm 1$ so the image of $M$ in $\PGL_2(\F_\el)$ is the identity. It follows now that $K_{k,\el}$ is totally real. 
One could arrive at a contradiction by approximating the roots of $P_{k,\el}$ to a high precision, 
but to get a proof one should use only symbolic calculations. 
The fields $K_{k,\el}$ with $\el\equiv 3\mod 4$ have negative discriminant hence cannot be totally real. Now suppose that a polynomial
$P(x)=x^n+a_{n-1}x^{n-1}+\cdots+a_0$ has only real roots. Then $a_{n-1}^2-2a_{n-2}$, being the sum of the squares of the roots, 
is non-negative and for a similar reason $a_1^2-2a_0a_2$ is non-negative as well. 
One can verify immediately that each of the polynomials $P_{k,\el}$ with $\el\equiv 1\mod 4$ fails at least one of these two criteria,   
hence none of the fields $K_{k,\el}$ involved in this paper is totally real. This proves the oddness of the representations $\rho_{k,\el}$.
Of course, this can also be checked with more general methods, like considering the trace pairing on $K_{k,\el}$ or invoking Sturm's theorem
\cite[Theorem 5.4]{JacAl}.
\newline

So now that we have verified all the conditions of Theorem \ref{serreconj} we remark as a final step that all spaces of modular forms $S_k(\Gamma(1))$ involved here are 
$1$-dimensional. So the modularity of each $\rho_{k,\el}$ implies immediately the isomorphism $\rho_{k,\el}\cong\rho_{\Delta_k,\el}$, 
hence also $\tilde\rho_{k,\el}\cong\tilde\rho_{\Delta_k,\el}$ , which completes the proof of Proposition \ref{PGLprop}.

\section{Proof of Corollary \ref{Lehmer}}\label{secLehmer}
If $\tau$ vanishes somewhere, then the smallest positive integer $n$ for which $\tau(n)$ is zero is a prime (see \cite[Theorem 2]{Le}). 
Using results on the exceptional representations for $\tau(p)$, Serre pointed out \cite[Section 3.3]{SeLa} that if $p$ is a prime number
with $\tau(p)=0$ then $p$ can be written as
$$
p=hM-1
$$
with 
$$
\eqalign{
&M=2^{14}3^75^3691=3094972416000,\cr
&\jac{h+1}{23}=1\quad\mbox{and}\quad h\equiv 0,30\mbox{ or }48\mod 49.
}
$$
In fact $p$ is of this form if and only if $\tau(p)\equiv 0\mod 23\cdot 49\cdot M$ holds.
Knowing this, we will do a computer search on these primes $p$ and verify whether \mbox{$\tau(p)\equiv 0\mod\el$} for $\el\in\{11,13,17,19\}$. 
To do this we need the following lemma.
\begin{lemma}
Let $K$ be a field of characteristic not equal to $2$. Then the following conditions on $M\in\GL_2(K)$ are equivalent:
\begin{enumerate}
\item[(1)]
$\tr M=0$.
\item[(2)]
For the action of $M$ on $\P^1(K)$, there are $0$ or $2$ orbits of length $1$ and all other orbits have length $2$.
\item[(3)]
The action of $M$ on $\P^1(K)$ has an orbit of length $2$.
\end{enumerate}
\end{lemma}
\begin{proof}
We begin with verifying (1)$\,\Rightarrow\,$(2). 
Suppose $\tr M=0$. Matrices of trace $0$ in $\GL_2(K)$ have distinct eigenvalues in $\o K$ because of $\tbox{char}(K)\not=2$. 
It follows that two such matrices are conjugate 
if and only if their characteristic polynomials coincide.
Hence $M$ and $M':= \mat{0}{1}{-\det M}{0}$ are conjugate so without loss of generality we assume $M=M'$.
Since $M^2$ is a scalar matrix, all the orbits of $M$ on $\P^1(K)$ have length $1$ or $2$. 
If there are at least $3$ orbits of length $1$ then $K^2$ itself is 
an eigenspace of $M$ hence $M$ is scalar, which is not the case. If there is exactly one orbit of length $1$ then $M$ has a non-scalar Jordan block
in its Jordan decomposition, which contradicts the fact that the eigenvalues are distinct.
\newline

The implication 
(2)$\,\Rightarrow\,$(3) is trivial so that leaves proving 
(3)$\,\Rightarrow\,$(1).
Suppose that $M$ has an orbit of order $2$ in $\P^1(K)$. After a suitable conjugation, we may assume that this orbit is 
$\{[{1 \choose 0}],[{0 \choose 1}]\}$. But this means that $M\sim\mat{0}{a}{b}{0}$ for certain $a,b\in K$ hence $\tr M=0$.  
\end{proof}
In view of this lemma it follows from Proposition \ref{PGLprop} that for $\el\in\{11,13,17,19\}$ and $p\not=\el$ we have that $\tau(p)\equiv 0\mod\el$ if and only if 
the prime $p$ decomposes in the number field $\Q[x]/(P_{12,\el})$ as a product of primes of degree $1$ and $2$, with degree $2$ occuring at least once.
For $p\nmid\disc(P_{12,\el})$, a property that all primes $p$ satisfying Serre's criteria possess, we can verify this condition by checking whether 
$P_{12,\el}$ has an irreducible factor of degree $2$ over $\F_p$; this can be easily checked by verifying 
$$
\o{x}^{p^2}=\o{x}\quad\mbox{and}\quad \o{x}^p\not= \o{x}
\quad\mbox{in}\quad \F_p[x]/(\o{P}_{12,\el}).
$$
Having done a computer search, it turns out that the first few primes satisfying Serre's criteria as well as $\tau(p)\equiv 0\mod 11\cdot 13\cdot 17\cdot 19$ are
$$
22798241520242687999,\ 60707199950936063999,\ 93433753964906495999.
$$
\section{The table of polynomials}
In this section we present the table of polynomials that is referred to throughout the article.
\begin{longtable}{|c|l|}
\caption {Polynomials belonging to projective modular representations}\label{polytable}\cr

\hline{}&{}\cr \multicolumn{1}{|c|}{$(k,\el)$} & \multicolumn{1}{c|}{$P_{k,\el}$}\cr{}&{}\cr \hline 
\endfirsthead

\multicolumn{2}{c}%
{{\tablename\ \thetable{} -- continued from previous page}} \cr
\hline{}&{}\cr \multicolumn{1}{|c|}{$(k,\el)$} &
\multicolumn{1}{c|}{$P_{k,\el}$}\cr{}&{}\cr\hline 
\endhead

\multicolumn{2}{r}{\footnotesize{Continued on next page}}
\endfoot

\hline
\endlastfoot

$(12,11)$&$x^{12} - 4x^{11} + 55x^9 - 165x^8 + 264x^7 - 341x^6 + 330x^5$\\*&$ -\, 165x^4 - 55x^3 + 99x^2 - 41x - 111$\cr
\hline
$(12,13)$&$x^{14} + 7x^{13} + 26x^{12} + 78x^{11} + 169x^{10} + 52x^9 - 702x^8 - 
1248x^7$ \\*& $+\,494x^6 + 2561x^5 + 312x^4 - 2223x^3 + 169x^2 + 506x - 215$\cr
\hline
$(12,17)$&$x^{18} - 9x^{17} + 51x^{16} - 238x^{15} + 884x^{14} - 2516x^{13} + 5355x^{12}$\\*&$ -\, 7225x^{11} - 1105x^{10} + 37468x^9 - 111469x^8 +
    192355x^7$\\*&$ -\, 211803x^6 + 134793x^5 - 17323x^4 - 50660x^3 + 47583x^2$\\*&$ -\, 19773x + 3707$
\cr\hline
$(12,19)$&$x^{20} - 7x^{19} + 76x^{17} - 38x^{16} - 380x^{15} + 114x^{14} + 1121x^{13}$\\*&$ -\, 798x^{12} - 1425x^{11} + 6517x^{10} + 152x^9 -
    19266x^8 - 11096x^7$\\*&$ +\, 16340x^6 + 37240x^5 + 30020x^4 - 17841x^3 - 47443x^2 $\\*&$-\, 31323x - 8055$\cr
\hline
$(16,17)$&$x^{18} - 2x^{17} - 17x^{15} + 204x^{14} - 1904x^{13} + 3655x^{12} + 5950x^{11}$\\*&$ -\, 3672x^{10} - 38794x^9 + 19465x^8 + 95982x^7 -
    280041x^6$\\*&$ -\, 206074x^5 + 455804x^4 + 946288x^3 - 1315239x^2 + 606768x $\\*&$-\, 378241$\cr
\hline
$(16,19)$&$x^{20} + x^{19} + 57x^{18} + 38x^{17} + 950x^{16} + 4389x^{15} + 20444x^{14}$\\*&$ +\, 84018x^{13} + 130359x^{12} - 4902x^{11} - 93252x^{10}
    + 75848x^9$\\*&$ -\, 1041219x^8 - 1219781x^7 + 3225611x^6 + 1074203x^5 $\\*&$-\, 3129300x^4 - 2826364x^3 + 2406692x^2 +
    6555150x - 5271039$\cr
\hline
$(16,23)$&$x^{24} + 9x^{23} + 46x^{22} + 115x^{21} - 138x^{20} - 1886x^{19} + 1058x^{18} $\\*&$+\, 59639x^{17} + 255599x^{16} +
    308798x^{15} - 1208328x^{14}$\\*&$ -\, 6156732x^{13} - 10740931x^{12} + 2669403x^{11} + 52203054x^{10}$\\*&$ +\, 106722024x^9
    + 60172945x^8 - 158103380x^7 - 397878081x^6$\\*&$ -\, 357303183x^5 + 41851168x^4 + 438371490x^3 +
    484510019x^2$\\*&$ +\, 252536071x + 55431347$\cr
\hline
$(18,17)$&$x^{18} - 7x^{17} + 17x^{16} + 17x^{15} - 935x^{14} + 799x^{13} + 9231x^{12}$\\*&$ -\, 41463x^{11} + 192780x^{10} + 291686x^9 -
    390014x^8 + 6132223x^7$\\*&$ -\, 3955645x^6 + 2916112x^5 + 45030739x^4 - 94452714x^3 $\\*&$+\, 184016925x^2 -
    141466230x + 113422599$\cr
\hline
$(18,19)$&$x^{20} + 10x^{19} + 57x^{18} + 228x^{17} - 361x^{16} - 3420x^{15} + 23446x^{14}$\\*&$ +\, 88749x^{13} - 333526x^{12} - 1138233x^{11} +
    1629212x^{10} $\\*&$+\, 13416014x^9 + 7667184x^8 - 208954438x^7 + 95548948x^6$\\*&$ +\, 593881632x^5 - 1508120801x^4 -
    1823516526x^3 $\\*&$+\, 2205335301x^2 + 1251488657x - 8632629109$\cr    
\hline
$(18,23)$&$x^{24} + 23x^{22} - 69x^{21} - 345x^{20} - 483x^{19} - 6739x^{18} + 18262x^{17}$\\*&$ +\, 96715x^{16} - 349853x^{15} +
    2196684x^{14} - 7507476x^{13}$\\*&$ +\, 59547x^{12} + 57434887x^{11} - 194471417x^{10} + 545807411x^9 $\\*&$+\, 596464566x^8
    - 9923877597x^7 + 33911401963x^6 $\\*&$-\, 92316759105x^5 + 157585411007x^4 - 171471034142x^3 
    $\\*&$+\, 237109280887x^2 - 93742087853x + 97228856961$\cr
\hline
$(20,19)$&$x^{20} - 5x^{19} + 76x^{18} - 247x^{17} + 1197x^{16} - 8474x^{15} + 15561x^{14}$\\*&$ -\, 112347x^{13} + 325793x^{12} -
    787322x^{11} + 3851661x^{10}$\\*&$ -\, 5756183x^9 + 20865344x^8 - 48001353x^7 + 45895165x^6 $\\*&$ -\, 245996344x^5 +
    8889264x^4 - 588303992x^3 - 54940704x^2 $\\*&$-\, 538817408x + 31141888$\cr
\hline
$(20,23)$&$x^{24} - x^{23} - 23x^{22} - 184x^{21} - 667x^{20} - 5543x^{19} - 22448x^{18}$\\*&$ +\, 96508x^{17} + 1855180x^{16} +
    13281488x^{15} + 66851616x^{14} $\\*&$+\, 282546237x^{13} + 1087723107x^{12} + 3479009049x^{11} $\\*&$+\, 8319918708x^{10} +
    8576048755x^9 - 19169464149x^8 $\\*&$-\, 111605931055x^7 - 227855922888x^6 - 193255204370x^5 $\\*&$+\, 
    176888550627x^4 + 1139040818642x^3 + 1055509532423x^2 $\\*&$ +\, 1500432519809x + 314072259618$\cr
\hline
$(22,23)$&$x^{24} - 11x^{23} + 46x^{22} - 1127x^{20} + 6555x^{19} - 7222x^{18}$\\*&$ -\, 140737x^{17} + 1170700x^{16} - 2490371x^{15} -
    16380692x^{14} $\\*&$+\, 99341324x^{13} + 109304533x^{12} - 2612466661x^{11} $\\*&$+\, 4265317961x^{10} + 48774919226x^9 -
    244688866763x^8$\\*&$ -\, 88695572727x^7 + 4199550444457x^6 - 10606348053144x^5 $\\*&$-\, 25203414653024x^4 +
    185843346182048x^3 $\\*&$-\, 228822955123883x^2 - 1021047515459130x $\\*&$+\, 2786655204876088$\cr

\end{longtable}

\end{document}